\documentclass{amsart}
\usepackage{amssymb,latexsym}
\usepackage{amsfonts}
\usepackage{amscd}
\usepackage{amsmath}
\usepackage{graphicx,array}
\usepackage[all]{xy}
\theoremstyle{plain}
\numberwithin{equation}{section}
\newtheorem{theorem}{Theorem}[section]

\newtheorem{lemma}[theorem]{Lemma}

\newtheorem{observation}[theorem]{Observation}
\theoremstyle{definition}

\newtheorem{remark}[theorem]{Remark}

\newcommand{\cA}{\mathcal{A}}

\newcommand{\cF}{\mathcal{F}}
\newcommand{\cO}{\mathcal{O}}
\newcommand{\cH}{\mathcal{H}}

\newcommand{\cK}{\mathcal{K}}
\newcommand{\cP}{\mathcal{P}}
\newcommand{\cR}{\mathcal{R}}

\newcommand{\conv}{\mathrm{conv}}

\newcommand{\fa}{\mathfrak{a}}

\newcommand{\fg}{\mathfrak{g}}

\newcommand{\fk}{\mathfrak{k}}

\newcommand{\fm}{\mathfrak{m}}
\newcommand{\fn}{\mathfrak{n}}

\newcommand{\fp}{\mathfrak{p}}

\newcommand{\fz}{\mathfrak{z}}
\newcommand{\gs}{\sigma}
\newcommand{\gD}{\Delta}

\newcommand{\R}{\mathbb{R}}
\newcommand{\C}{\mathbb{C}}

\newcommand{\Tr}{\mathrm{Tr} }

\newcommand{\id}{\mathrm{id}}

\renewcommand{\Re}{\mathrm{Re}}
\renewcommand{\Im}{\mathrm{Im}}

\newcommand{\Cr}{\mathrm{Cr}(X)}

\newcommand{\XiO}{\Xi (\Omega )}

\theoremstyle{plain}
\newcommand{\ind}{\mathrm{ind}}
\newcommand{\cHC}{\cH_X}
\newcommand{\cHXi}{\cH_{\Xi}}

\begin{document}
\title[Representation theory, Radon transform and the heat equation]
{Representation theory, Radon transform and the heat equation
on a Riemannian symmetric space}

\author{Gestur \'{O}lafsson and Henrik Schlichtkrull}

\address{Department of Mathematics, Louisiana State
University, Baton Rouge,
LA 70803, USA} \email{olafsson@math.lsu.edu}
\address{Department of Mathematics, University of Copenhagen,
Universitetsparken 5, DK-2100 K{\o}benhavn {\O}, Denmark}
\email{schlicht@math.ku.dk}
\thanks{Research of \'Olafsson was supported by NSF grants
  DMS-0402068.}
\thanks{Part of this paper was written while the second author was
a visitor at the Max-Planck Institute for Mathematics, Bonn, Germany.
He expresses gratitude for support and hospitality}

\begin{abstract}
Let $X=G/K$ be a Riemannian symmetric space of the noncompact type.
We give a short exposition of the representation theory
related to $X$, and discuss its holomorphic extension to the complex crown, a
$G$-invariant subdomain in the complexified
symmetric space $X_\C=G_\C/K_\C$. Applications to the heat
transform and the Radon transform for $X$ are given.
\end{abstract}

\maketitle
\section*{Introduction}
\noindent
In the analysis of Riemannian symmetric spaces
one can follow several approaches, with the emphasis
for example
on differential geometry, partial differential equations,
functional analysis, complex analysis, or representation theory.
The representation theory associated with a
Riemannian symmetric space
is of course well known through the work of Harish-Chandra
\cite{HC57,HC58}, and it is based on important results in the
theory of infinite dimensional representations, by Gelfand-Naimark,
Mackey, Segal, just to mention a few.

The powerful theory of Harish-Chandra has been widely generalized
\cite{BS05,D98,HC76},
and successfully applied to many problems
related to Riemannian symmetric spaces. The geometric point of
view, with emphasis on relations to topics like classical
Euclidean analysis and the Radon transform has been
represented by Helgason \cite{He63,He68,He70}.
Both aspects, as well as the connection
to the work of the school around Gelfand and Graev are
described in a short and clear fashion in the second half of
Mackey's famous Chicago Lectures on unitary group representations
\cite{M76}.

In recent years much research has been directed towards the
interplay between the real analysis and the geometry of the symmetric space
$G/K$ on the one side, and complex methods in analysis,
geometry and representation
theory on the other side. This development has a long history tracing back to
Cartan's analysis of bounded symmetric domains, and Harish-Chandra's
construction of the holomorphic discrete series \cite{HC55}.
The so-called Gelfand-Gindikin program \cite{GG77},
suggests to consider functions on $G$ through
holomorphic extension to domains in the complexification $G_{\mathbb C}$
of $G$, and to study representations of $G$ realized
on spaces of such functions, analogous to the classical Hardy spaces
on tube domains over $\R^n$. Only partial results have been obtained
so far. The program was carried out for the holomorphic discrete series
of groups of Hermitian type type in  \cite{GIO82a,RS86},
the holomorphic discrete series for compactly causal symmetric spaces
in \cite{HOO91}, and finally the holomorphic most continuous
series for noncompactly
causal symmetric spaces in \cite{GKO04a}.

Connected to this development
is the study of natural domains
in $G_{\mathbb C}$ to which the spherical function of $G/K$
admit holomorphic extensions, initiated by
Akhiezer and  Gindikin \cite{AG90}, and continued by
several people.
Without claiming to be complete
we would like to mention
\cite{B03,BHH03,F03a,JF03,FH02,FHW06,GKO03,GKO04a,GKO06a,GM01,KS04,KS05,M06}
representing different aspects of this important development.
The most relevant articles related to
the present exposition are the papers \cite{GK02,KS05}.
One important conclusion is that there exists a maximal
$G$-invariant domain, called the {\it complex crown},
to which all the spherical functions extend.

A related area of recent research is the study of the
heat equation on the symmetric space.
There are several generalizations to the well-known Segal-Bargmann transform,
which maps an $L^2$-function $f$ on the Euclidean space $\R^n$ to
the function $H_tf$ on $\C^n$, which is the holomorphic continuation of
the solution at time $t$ to the heat equation with initial
Cauchy data $f$.
The first work in this direction was by Hall \cite{H94}, who replaced
$\R^n$ by a connected compact semisimple Lie group $U$, and
$\C^r$ by the complexification $U_\C$.
This was put into a
general framework using polarization of a restriction map in \cite{OO96}.
The results of Hall were
extended to compact
symmetric spaces $U/K$
by Stenzel in \cite{S99}, where the complexification is  $U_\C/K_\C$.
It is important to note that
in this compact case, all eigenfunctions of the algebra of invariant
differential operators on $U/K$, and also the heat kernel itself,
extend to holomorphic functions on $U_\C/K_\C$.
This is related to the fact, that each irreducible representation
of $U$ extends to a holomorphic representation of $U_\C$.

For a symmetric space of the noncompact type  $G/K$,
the maximal $G$-invariant domain is the complex crown and
not the full complexification $X_\C=G_\C/K_\C$.
It was  shown
in \cite{KOS05} that the image of the Segal-Bargmann transform on
$G/K$ can be identified as a Hilbert-space of holomorphic functions on
the crown. The norm was defined using orbital integrals and
the Faraut version of Gutzmer's formula \cite{JF03}. Some special
cases have also been considered in \cite{H04a,H04b},
but without using the Akhiezer-Gindikin domain explicitly. In
particular, in \cite{H04b} Hall and Mitchell
give a description of
the image of the $K$-invariant functions in $L^2(G/K)$ in case
$G$ is complex.

In the present paper we have the following basic aims.
Our first aim is to give
a short exposition of the basic representation theory related
to a Riemannian symmetric space $G/K$, connecting to work of
Harish-Chandra and Helgason. The secondary aim is to
discuss the complex extension in $X_\C$,
and to
introduce a $G$-invariant Hilbert space of holomorphic functions
on the crown,
which carries all the representation theory of $G/K$.
Thirdly, we combine the Fourier theory and the holomorphic
theory in a study of the heat equation on $G/K$.
The image of the Segal-Bargmann transform is described
in a way similar to that of  \cite{OS06},
where we considered the case of $K$-invariant functions
on $G/K$ (in a more general setting of root systems),
and derived a result containing that of
Hall and Mitchell as a special case.
The image is characterized
in terms of a Fock space on $\fa_\C$.

Finally we discuss some aspects of
the Radon transform on $G/K$.
We introduce a $G$-invariant Hilbert space of CR-functions
on a subset of the complexified horocycle space,
and we use the Radon transform to construct
a unitary $G$-invariant isomorphism
between the Hilbert spaces
on the crown and on the complex horocycle space, respectively.

We shall now describe the content of the paper in some more detail.
Let $X=G/K$ be a Riemannian symmetric space of the noncompact type.
For the mentioned primary aim, we define a Fourier transform and we
state a Plancherel theorem and an inversion formula, following the
beautiful formulation of Helgason. The content of the formulation
is described by means of representation theory.
We also indicate
a representation theoretic proof of the reduction to
Harish-Chandra's theorem.
At the end of Section \ref{se-not} we give a short description of
the heat transform and its image in $L^2(X)$ as a direct integral
over principal series representations.

The space $X$  is naturally contained in the complexified symmetric
space $X_{\mathbb C}=G_{\mathbb C}/K_{\mathbb C}$
as a totally real submanifold.
Inside the complexification $X_{\mathbb C}$ is the $G$-invariant
domain $\Cr\subset X_{\mathbb C}$, the {\it complex crown},
which was introduced in \cite{AG90}. It has
the important property, shown in \cite{KS05},
that all joint eigenfunctions for $\mathbb D(X)$,
the algebra of invariant differential operators on $X$,
extend to holomorphic
functions on it. A second fundamental fact is the
convexity theorem of Gindikin, Kr\"otz, and Otto \cite{GK02,KOt06},
which is recalled in Theorem \ref{th-convexity}.
In Section \ref{s-cHC}
we define a $G$-invariant Hilbert space
$\cHC$ of holomorphic functions on $\Cr$, such that
restriction to $X$ maps continuously into $L^2(X)$, and such
that the representation of $G$ on  $\cHC$ carries
all the irreducible representations found in the decomposition
of $L^2(X)$. The definition, which is essentially
representation theoretic, is related to definitions given in
\cite{GKO03}. The main results are stated in Theorem \ref{th-3.1}
where a representational description of $\cHC$ is given and
the reproducing kernel of the space is determined, see also
\cite{GKO06a}.

The final section deals with the holomorphic extension of $H_tf$ for
$f\in L^2(X)$.
It was shown in \cite{KS05} that $H_tf$ extends to a holomorphic function
on $\Cr$, and a description of the image
$H_t(L^2(X))\subset \cO (\Cr)$ was given in
\cite{KOS05}. In fact, we show in Theorem \ref{th-ImageHt2} that
$H_t(L^2(X))\subset \cHC$, and we give an alternative description of
the image by means of the Fourier transform.

The Radon transform sets up a relation
between functions on $X$ and functions on the space
$\Xi=G/MN$ of horocycles on $X$. The
Radon transform
is used in Theorem \ref{th-image2} to give yet another
description of the image of $H_t$. We  also give an inversion
formula for $H_t$.

The space $\Xi$ sits inside a complex space, the space $\Xi_{\mathbb C}$
of complex horospheres on $X_{\mathbb C}$, and in this complex space
one can define a $G$-invariant domain $\XiO\supset \Xi$ analogous
to the crown. However, $\XiO$ is not a complex manifold but a CR-submanifold
of $\Xi_{\mathbb C}$. On this domain, we define a space $\cHXi$
similar to $\cHC$, and we show in Theorem \ref{th-HCIsoHXI}
that the normalized
Radon transform can be used to set up an
unitary isomorphism $\tilde{\Lambda} : \cHC\to \cHXi$.

\section{Representation theory and harmonic analysis }\label{se-not}
\noindent
In this section we introduce the standard notation
that will be used throughout this article. We will also
recall some well known facts about the principal series
representations $\pi_\lambda$ and the Fourier transform on $X$.

\subsection{Notation}\label{ss-str}
Let $G$ be a connected noncompact semisimple Lie group with
Lie algebra $\fg$, let $\theta :G\to G$ be a Cartan involution and set
$$K=G^\theta=\{x\in G\mid \theta (x)=x\}\, .$$
The space $X=G/K$ is Riemannian symmetric space of noncompact type.
It does not depend on which one of the locally isomorphic groups
with Lie algebra $\fg$ is used, as the center of $G$ is always contained
in $K$. We will therefore assume that $G$ is contained in
a simply connected Lie group $G_\C$ with Lie algebra $\fg_\C=
\fg\otimes_\R \C$. In particular,  $G$ has finite center and
$K$ is a maximal compact subgroup of $G$.

We will denote by the same symbol $\theta$
the holomorphic extension of $\theta$ to $G_\C$, as well as
the derived Lie algebra homomorphisms $\theta: \fg\to \fg$
and its complex linear extension to $\fg_\C$.
Denote by $\sigma : \fg_\C\to \fg_\C$ the conjugation
on $\fg_\C$ with respect to $\fg$.
As $G_\C$ is simply connected, $\gs$ integrates to
an involution on $G_\C$ with $G=G_\C^\sigma$.
Denote by $K_\C\subset G_\C$ the complexification
of $K$ in $G_\C$ and $X_\C =G_\C/K_\C$.
Then $X_\C$ has a complex structure, with respect to which
$X\simeq G\cdot x_o\subset X_\C$ is a totally real
submanifold. It can be realized as the connected component
containing  $x_o=eK_\C\in X_\C$ of the fixed point set of
the conjugation $g\cdot x_o\mapsto \sigma (g\cdot x_o)
:=\sigma (g)\cdot x_o$.

Let
$\fg=\fk\oplus \fp$ be the Cartan decomposition
defined by $\theta$, and let $\fa$ be a maximal abelian
subspace of $\fp$, $\Delta$ the set of roots of
$\fa$ in $\fg$, and $\Delta^+$ a fixed set of positive
roots. For $\alpha\in\Delta$ let
$$\fg^\alpha=\{X\in \fg\mid (\forall H\in \fa)\, [H,X]=\alpha (H)X\}$$
be the joint $\alpha$-eigenspace. Let
$$\fn=\bigoplus_{\alpha\in\Delta^+}\fg^\alpha\, ,
$$
then $\fn$ is a nilpotent subalgebra of $\fg$.
Let
$\fm =\fz_\fk (\fa)$,
then
$\fm\oplus \fa\oplus\fn $ is a minimal parabolic subalgebra and
$$\fg=\fk\oplus \fa\oplus\fn
\, .$$
Let $A=\exp\fa$, $N=\exp\fn$, $M=Z_K(A)$ and $P=MAN$. We have the Iwasawa decompositions
$$G=KAN \subset K_\C A_\C N_\C\subset G_\C$$
where the subscript ${}_\C$ on a subgroup of $G$ stands
for its complexification in $G_\C$.
The map $K\times A\times N\ni (k,a,n)\mapsto kan \in G$ is
an analytic diffeomorphism (the analogous statement
fails for the complexified Iwasawa decomposition).
We shall denote its inverse
$x\mapsto (k(x),a(x),n(x))$.

Let $B=K/M$, this is the
so-called Furstenberg boundary of $X$.
Note that since $B\simeq G/P$, it carries
the action of $G$ given by $g\cdot (kM)=k(gk)M$.

\subsection{Integration}\label{ss-int}
If $C$ is a Lie group or a homogeneous space
of a Lie group, then we
denote by
$dc$ a left invariant measure on
$C$. We normalize the invariant measure on compact groups
and compact homogeneous spaces
such that the total measure is one. We require
that
the Haar measures on $A$ and $i\fa^*$ are normalized such that
if
\begin{equation}\label{eq-FTonA}
\cF_A (f)(\lambda )=\int_A f(a)a^{-\lambda} \, da
\end{equation}
is the Fourier transform on the Abelian group $A$, then
$$f(x)=\int_{i\fa^*} \cF_A (f)(\lambda )a^{\lambda}\, d\lambda\, .$$
Finally, we normalize
the Haar measure on $N$ as usual
by $\int_N a(\theta (n))^{-2\rho }\, dn=1$,
where  $2\rho = \sum_{\alpha \in \gD^+}m_\alpha \alpha$ and
$m_\alpha =\dim \fg^\alpha$. Then
we can normalize the Haar measure on $G$ such that for all $f\in C_c(G)$ we
have
\begin{equation}\label{eq-HaarOnG}
\int_G f(x)\, dx=\int_N\int_A\int_K f(ank)\, dkdadn
=\int_K\int_A\int_N f(kan)a^{2\rho}\, dndadk\, .
\end{equation}
In this normalization the invariant measure on $X$ is given by
$$\int_X f(x)\, dx = \int_G f(g\cdot x_o)\, dg =\int_N\int_A f(an\cdot x_o)\, dadn\, .$$
Finally, it follows from (\ref{eq-HaarOnG}) that the $K$-invariant
measure on $B$ transforms under the $G$-action according to
\begin{equation}\label{eq-GactionOnB}
\int_B f(g\cdot b)\, db=\int_B f(b)a (g^{-1}b)^{-2\rho}\, db
\end{equation}
for all $f\in L^1(B)$ and $g\in G$.

\subsection{Spherical principal series  and spherical functions}\label{ss-sp}
Denote by $L$ the left regular representation
$L_af(x)=f(a^{-1}x)$ and by $R$ the right regular
representation $R_af(x)=f(xa)$. If $C/D$ is a homogeneous space,
then we identify functions on $C/D$ with right $D$-invariant functions
on $C$.

The {\it spherical principal series} of representations
are defined as follows.
For $\lambda\in \fa_\C^*$ denote by $H_\lambda$
the
Hilbert space of measurable functions $f: G\to \C$ such that
for all  $man\in P$ (with the obvious notation)
\begin{equation}\label{eq-ind}
R_{man} f = a^{-\lambda -\rho}f\quad \mathrm{and}
\quad \|f\|^2=\int_K |f(k)|^2\, dk<\infty\, .
\end{equation}
Define the representation $\pi_\lambda$ of $G$ on $H_\lambda$ by
\begin{equation}\label{eq-pi}
\pi_\lambda (x)f(y)= L_xf(y) = a(x^{-1}y)^{-\lambda - \rho}f(k(x^{-1}y))\, .
\end{equation}

A different picture (called the \textit{compact picture})
of the representations $\pi_\lambda$ is obtained by noting that
the Iwasawa decomposition implies that the restriction map
$f\mapsto f|_K$,
is a unitary $K$-isomorphism $H_\lambda \to L^2(B)$.
Thus we may view $\pi_\lambda$ as
a representation on the latter space.
By (\ref{eq-pi}) the representation is then given by
\begin{equation}\label{eq-compactpicture}
\pi_\lambda (x)f(b)=a(x^{-1}k)^{-\lambda - \rho}f(x^{-1}\cdot b)\, .
\end{equation}
The advantage of the compact picture is that
the Hilbert space is independent
of $\lambda$.

The representations $\pi_\lambda$ are known to be unitary
when $\lambda$ is purely imaginary on $\fa$. Furthermore,
it is known that $\pi_\lambda$ is irreducible
for almost all
$\lambda\in \fa_\C^*$, and that $\pi_{w\lambda}$ and
$\pi_\lambda$ are equivalent representations, also
for almost all $\lambda\in\fa_\C^*$, see \cite{B56,K69}.
Here $W$ denotes the
Weyl group of the root system $\Delta$.

The representations $\pi_\lambda$ all restrict to the same representation
of $K$, the left regular representation on $L^2(B)$.
In particular, the trivial representation of $K$ has multiplicity
one and is realized on the space of constant functions
on $B$. We fix $p_\lambda\in H_\lambda^K$ as the
constant function 1 on $B$, that is, on $G$ it is
$$p_\lambda (g)=a (g)^{-\lambda-\rho}.$$
We define the following function $e_{\lambda,b}$
on $X$ for
$(\lambda,b)\in \fa_\C^*\times B$,
\begin{equation}\label{eq-lambda}
e_{\lambda, kM}(gK)=p_\lambda (g^{-1}k)
\, ,\quad
k\in K, g\in G\, .
\end{equation}
then $(x,b)\mapsto e_{\lambda,b}(x)$ is the generalized {\it Poisson kernel}
on $X\times B$.
The {\it spherical functions} on $X$ are the $K$-biinvariant
matrix
coefficients of $\pi_\lambda$ defined by
\begin{equation}\label{spherical function}
\varphi_{\lambda }(x)=
(\pi_\lambda (g)p_\lambda ,p_\lambda )
=\int_B e_{\lambda ,b}(x) \, db
\end{equation}
where $x=gK$ and
$\lambda\in\fa^*_\C$. The latter integral is exactly
Harish-Chandra's formula for the spherical functions.

\subsection{The standard intertwining operators}\label{ss-intop}
As mentioned, the representation $\pi_{w\lambda}$
is known to be equivalent with $\pi_\lambda$ for almost all
$\lambda\in i\fa^*$. Hence for such $\lambda$ there exists a
unitary intertwining operator
$$\cA(w,\lambda)\colon H_\lambda \to H_{w\lambda}.$$
The operator is unique, up to scalar multiples,
by Schur's lemma. The trivial $K$-type
has multiplicity one and is generated by the function $p_\lambda\in H_\lambda$,
which has norm 1 in $L^2(B)$.
It follows that $\cA(w,\lambda)p_\lambda$ is a unitary multiple
of $p_{w\lambda}$. We normalize the intertwining operator  so that
$$\cA(w,\lambda)p_\lambda=p_{w\lambda}.$$
The intertwining operator so defined is called the
{\it normalized standard intertwining operator}.
It is known that the map $\lambda\mapsto\cA(w,\lambda)$
extends to a rational map (which we denote by the same symbol)
from $\mathfrak a^*_{\mathbb C}$ into the bounded operators on $L^2(B)$.

In fact, one can give a formula for the operator
$\cA(w,\lambda)$ as follows.
The (unnormalized) {\it standard intertwining operator}
$A(w,\lambda)\colon H_\lambda \to H_{w\lambda}$
is defined by the formula
$$A(w,\lambda)f(g)=\int_{\bar N_w} f(gw\bar n)
\,d\bar n $$
where $\bar N_w=\theta( N)\cap w^{-1}Nw$, see \cite{KS71}.
The integral converges when $f$ is continuous and
$\lambda\in\mathfrak a^*_{\mathbb C}$ satisfies
that ${\rm Re}\langle\lambda,\alpha\rangle> 0$
for all $\alpha\in\Delta^+$. It is defined by
meromorphic continuation
for other values of $\lambda$, and by continuous extension
for $f\in L^2(B)$ (see \cite{Knapp}, Ch. 7).
{}From the definition of $A(w,\lambda)$ we see
$$A(w,\lambda)p_\lambda= c_w(\lambda)p_{w\lambda}$$
where
$c_w(\lambda)=\int_{\bar N_w} p_\lambda(\bar n)
\,d\bar n $
(see \cite{He84} page 446). Hence
$$\cA(w,\lambda)=c_w(\lambda)^{-1}A(w,\lambda).$$

\subsection{The Fourier transform}\label{ss-ft}
In this section we introduce the Fourier transform on $X$, following
Helgason \cite{He68,He70}.
While Helgason introduced it from a more
geometric point of view, we shall show here that it can also be done
from a representation theory point of view, resulting in alternative
proofs of the inversion formula and the Plancherel theorem.

{}From the point of view of representation
theory, the Fourier transform of an integrable function on $G$
is the operator $\pi(f)$ on $\cH$ defined by
$$\pi(f)v=\int_G f(g)\pi(g)v\,dg, \quad v\in\cH,$$
for each unitary irreducible representation $(\pi,\cH)$.

If $f$ is a function on $X$, this operator will be trivial
on the orthocomplement of the space of
$K$-fixed vectors. This space
is always one dimensional in an irreducible representation,
and hence it becomes natural to define the Fourier transform of
$f$ as the vector $\pi(f)v$ in the representation space of $\pi$,
where $v$ is a specified $K$-fixed vector. For the spherical
principal series, we thus arrive at the following definition
of the Fourier transform:
$$\hat f_\lambda:=\pi_{-\lambda}(f)p_{-\lambda}\in H_{-\lambda}$$
for $\lambda\in i\fa^*$
(the reason for the minus is just historical).

In the notation of the compact picture, it is
\begin{equation}
\hat{f}_\lambda(b)=(\pi_{-\lambda} (f)p_{-\lambda} ) (b)
=\int_G f(x)p_{-\lambda} (x^{-1}b)\, dx,
=\int_G f(x)\overline{e_{\lambda ,b} (x)} \, dx\, .
\end{equation}
where $b\in B$.
Thus the Fourier transform of
$f$ may be viewed as a map
$$
\fa_\C^* \times B \ni (\lambda , b )\mapsto
\hat{f}( \lambda, b ):=\hat{f}_\lambda(b)\in \C$$
Apart from the replacement of $\lambda$ by $i\lambda$,
this is the Fourier transform as it was introduced by
Helgason in \cite{He68} .

Let
$$\fa^+:=\{H\in \fa\mid (\forall \alpha \in\Delta^+)\, \alpha (X)>0\}$$
be the
positive Weyl chamber corresponding to $\Delta^+$, and let
$\fa_+^*$ denote the corresponding open chamber in $\fa^*$.
Let $c(\lambda )$ be the Harish-Chandra $c$-function,
which for $\Re\lambda \in\fa^*_+$ is given by (see \cite{He84} p. 447)
$$c(\lambda)=c_{w^*}(\lambda)=\int_{\bar N} p_\lambda(\bar n) d\bar n,$$
where $w^*\in W$ is the long element and
$\bar N=\theta N$. We recall that an explicit formula for
$c(\lambda)$ was determined by Gindikin and Karpelevich,
see \cite{GK62} or \cite{He84}, p. 447.
Furthermore, we define a measure
$\mu$ on $i\fa^*_+\times B$
by
\begin{equation}\label{eq-mu}
d\mu ( \lambda ,b )=|c( \lambda )|^{-2}\, d\lambda \,db\, .
\end{equation}
We will also denote by $d\mu$ the measure $|c( \lambda
)|^{-2}d\lambda$ on $i\fa^*$.

Let $L^2_W( i\fa^*\times B,\frac{d\mu }{|W|} )$ be the space of
all $F\in L^2( i\fa^*\times B , \frac{d\mu }{|W|})$ such that
for all $w\in W$  we have
\begin{equation}\label{intertwining equation}
F(w\lambda,\cdot)=\cA(w,-\lambda)F(\lambda,\cdot)
\end{equation}
in $L^2(B)$, for almost all $\lambda \in i\fa^*$.
Notice that, since $\cA(w,-\lambda)$ is an intertwining operator
for each $\lambda\in i\fa^*$, this is an invariant subspace for
the unitary action of $G$ on
$L^2(i\fa^*\times B , \frac{d\mu }{|W|})$, defined
 by
$$(g\cdot F)(\lambda,\cdot)=\pi_{-\lambda}(g)F(\lambda,\cdot).$$

We recall the following theorem of \cite{He70,He73}:
\begin{theorem}\label{th-He}
The Fourier transform is an intertwining unitary  isomorphism
\begin{equation}\label{eq-isomorphism}
L^2(X)\simeq L^2_W(i\fa^*\times B ,\frac{d\mu }{|W|} )\, .
\end{equation}
Furthermore, if $f\in C_c^\infty (X)$, then
\begin{equation}\label{eq-inversion}
f(x)=\frac{1}{|W|}\int_{i\fa^*}\int_B \widehat{f}
(\lambda  ,b )e_{\lambda, b}(x)\, d\mu (\lambda  ,b)\, .
\end{equation}
\end{theorem}

\bigskip
For left $K$-invariant functions on $X$, this is
Harish-Chandra's Plancherel theorem for $X$. The Fourier transform
$\hat f_\lambda$ is then a constant function
on $B$, and the constant is $$\hat f_\lambda(b)=
\int_K \hat f(\lambda,kM)\,dk
=\int_X f(x)\varphi_{-\lambda}(x)\,dx$$
(see (\ref{spherical function})),  the
{\it spherical Fourier transform} of $f$. By definition,
the normalized intertwining operator maps $p_{-\lambda}$
to $p_{-w\lambda}$, and hence in this case the intertwining relation
(\ref{intertwining equation}) is just $F(w\lambda)=F(\lambda)$.

The proof of the spherical Plancherel theorem is given in
\cite{HC58} (see also \cite{HC54}),
but it depends on two conjectures, see p. 611-612.
One conjecture was affirmed with \cite{GK62} by the mentioned formula
for $c(\lambda)$. The second conjecture is affirmed
in \cite{HC66}, see p. 4. A simpler proof has
later been given in \cite{Ro77}, see also \cite{He84} p. 545.

As explained in \cite{He66} p. 50, Theorem \ref{th-He} is proved by
reduction to the spherical case.
Because of the modified point of view invoking the representation theory,
and since we have stated the relations (\ref{intertwining equation})
differently, we discuss some aspects of the proof.
For more details, we refer to \cite{He70,He73,He94}.

\subsection{The intertwining relation}\label{ss-intertwiningRel}
\noindent
For the Fourier transformed function $F=\hat f$,
the intertwining relation (\ref{intertwining equation})
is a direct consequence of the
definition of $\hat f$:
\begin{eqnarray}
\cA(w,-\lambda)\hat f_\lambda&=&
\cA(w,-\lambda)\pi_{-\lambda}(f)p_{-\lambda}\nonumber\\
&=&\pi_{-w\lambda}(f)\cA(w,-\lambda)p_{-\lambda}\nonumber\\
&=& \pi_{-w\lambda}(f)p_{-w\lambda}\nonumber\\
&=&
\hat f_{w\lambda}\label{eq-IntertwRel}\, .
\end{eqnarray}

The equations (\ref{intertwining equation})
allow the following important reformulation
(which is the original formulation of Helgason, see \cite{He73}, p. 132)
\begin{equation}\label{eq-invariance}
\int_B e_{w \lambda ,b } (x )F(w \lambda ,b )\, db
=\int_B e_{\lambda ,b }(x )F(\lambda )\, db\,
\qquad (\forall x\in X).
\end{equation}
It follows that the integral over $i\fa^*$ in (\ref{eq-inversion})
is $W$-invariant, and thus can be written as an integral over
the Weyl chamber $\fa^*_+$.

The equivalence of (\ref{intertwining equation})
and (\ref{eq-invariance}) is an immediate
consequence of the following lemma.
\begin{lemma}
Let $f,g\in L^2(B)$ and $w\in W$ be given.
Then the following holds for every $\lambda\in\mathfrak a^*_{\mathbb C}$
outside a locally finite set of complex hyperplanes.
The relation
\begin{equation}\label{intertwining relation}
\int_B e_{w\lambda ,b } (x )g(b )\, db
=\int_B e_{\lambda ,b }(x )f(b)\, db\, .
\end{equation}
holds for all $x\in X$,
if and only if
\begin{equation}\label{g=Af}
g=\cA(w,-\lambda)f.
\end{equation}
\end{lemma}

\begin{proof} The relation (\ref{intertwining relation})
can be written in terms of
the {\it Poisson transformation} for $X$.
Recall (see \cite{He70,Sch84})
that the Poisson transform is the operator
$\cP_\lambda\colon H_{-\lambda}\to C^\infty(X)$
which is defined by
$$\cP_\lambda f(x)=\int_B f(b) e_{\lambda,b}(x)\,db$$
or equivalently (see \cite{Sch84} page 80-81),
$$\cP_\lambda f(gK)=\int_K f(gk)\,dk.$$
The latter formula shows that $\cP_\lambda$ is a $G$-equivariant
operator for the left action.
Notice that $\cP_\lambda p_{-\lambda}$ is exactly
the spherical function $\varphi_\lambda$.
It is now seen that
(\ref{intertwining relation}) holds for all $x$ if and only if
\begin{equation}\label{W on Poisson}
\cP_{w\lambda}g=\cP_\lambda f.
\end{equation}

On the other hand, the following can be seen to hold
for all $\lambda$, for which the normalized
standard intertwining operator is non-singular,
\begin{equation}\label{intertwining Poisson}
\cP_{w\lambda}\circ\cA(w,-\lambda)= \cP_\lambda.
\end{equation}
Indeed, it suffices to verify the identity for
almost all $\lambda$, so we can assume that
$\pi_\lambda$ is irreducible. It follows from the identity
$\varphi_{w\lambda}=\varphi_\lambda$
that the two operators agree
when applied to the element $p_{-\lambda}\in H_{-\lambda}$.
Since the operators are $G$-equivariant, they must
then agree everywhere.
This proves (\ref{intertwining Poisson}).

The equivalence of (\ref{W on Poisson}) and (\ref{g=Af}) follows immediately,
for all $\lambda$ such that $\cP_{w\lambda}$ is injective. The latter
is obviously true if $\pi_{-w\lambda}$ is irreducible. In fact
it is known that the Poisson transformation is injective
for all $\lambda$, except on a singular set of hyperplanes
(see \cite{He76} and \cite{Sch84} Thm 5.4.3).
\end{proof}

\subsection{The inversion formula}
Let $f\in C^\infty_c(X)$. Then, at the origin of $X$,  $f(eK)=f^K(eK)$ where
$f^K(x)=\int_K f(kx)\, dk.$
By the inversion formula of Harish-Chandra one can determine
$f(eK)$  through the expression
\begin{equation}\label{eq-f(eK)}
f(eK)=\int_{i\fa_+^*} \Tr(\pi_{-\lambda }(f))
\, d\mu (\lambda).
\end{equation}
Applied to the function $L_{g^{-1}}f$ it gives
$$f(gK)=\int_{i\fa_+^*} \Tr (\pi_{-\lambda}(g^{-1})\pi_{-\lambda }(f))
\, d\mu (\lambda).$$
Since $\pi_{-\lambda }(f)$ annihilates the orthocomplement
of $H_{-\lambda}^K$,
$$\Tr (\pi_{-\lambda}(g^{-1})\pi_{-\lambda }(f))
=(\pi_{-\lambda }(g^{-1})\pi_{-\lambda }(f)p_{-\lambda}
,p_{-\lambda})_{L^2(B)},
$$
and  as the representations are unitary this equals
\begin{equation*}
( \pi_{-\lambda }(f)p_{-\lambda},\pi_{-\lambda }(g) p_{-\lambda })_{L^2(B)}
=\int_B \hat{f}(\lambda ,b)e_{\lambda ,b}(gK)\, db.
\end{equation*}
It follows that
$$f(x)=\int_{i\fa_+^*\times B}\hat{f}(\lambda ,b)e_{\lambda ,b}(x)\, d\mu (\lambda ,b)\, .$$

\subsection{The $L^2$-isomorphism}
For $f:G\to \C$ let $f^*(g)=\overline{f(g^{-1})}$. Then, if $\pi $ is an unitary
representation, we have $\pi (f^*)=\pi (f)^*$ and $\pi (f*g)=\pi (f)\pi (g)$ for
$f,g\in L^1(G)$.
Thus, for $\lambda\in i\fa^*$ (so $\pi_{\lambda}$ is unitary), we have
\begin{equation}\label{trace}
\Tr\, ( \pi_{-\lambda} (f^* *f))
=(\pi_{- \lambda} (f)p_{-\lambda} ,\pi_{-\lambda} (f)
 p_{-\lambda } )=\int_B |\hat{f} ( \lambda, b )|^2\, db.
\end{equation}
As
$\|f\|^2=f^* *f(eK)$ the equations (\ref{eq-f(eK)}) and (\ref{trace}) imply that
$$\|f\|_{L^2(X)}^2=\int_{i\fa^*_+ \times  B}|\hat{f}(\lambda,b )|^2\, d\mu (\lambda ,b )=
\int_{i\fa_+^*} \|\hat{f}_\lambda  \|^2\, d\mu (\lambda )\, .$$
Thus $L^2(X)\ni f\mapsto \hat{f}_\lambda \in L^2(B)$ sets up a
unitary map into
$\int_{i\fa_+^*}^\oplus(H_{-\lambda},\pi_{-\lambda})\, d\mu (\lambda )$.
By Harish-Chandra's Plancherel formula
this is an unitary isomorphism
on the level of $K$-invariant elements.
We want to show that this implies it is onto.
Let $\cK$ be the orthogonal complement of the image in
$\int_{i\fa_+^*}^\oplus H_{-\lambda}  d\mu (\lambda )\, .$
Then $\cK$ is $G$-invariant and $\cK^K=\{0\}$.
By \cite{M76}, Thm. 2.15,
there exists measurable subset $\Lambda\subset i\fa_+^*$  such that,
as a representation of $G$, we have
 $\cK\simeq \int_{\Lambda}^\oplus H_{-\lambda}\, d\mu (\lambda )$.
In particular,
 $\cK^K\simeq \int_{\Lambda}^\oplus H_{-\lambda}^K\, d\mu (\lambda )$
and hence  it follows from $\cK^K=\{0\}$
that $\mu(\Lambda)=0$. Hence $\cK=\{0\}$ and
\begin{equation}\label{eq-L2XasDirectIntegral}
(L^2(X),L)\simeq \int_{i\fa_+^*}^\oplus(H_{-\lambda},\pi_{-\lambda})\, d\mu (\lambda )
\end{equation}
The isomorphism statement in (\ref{eq-isomorphism}) follows
easily.
This completes the discussion of  Theorem \ref{th-He}.

\subsection{The heat equation}\label{ss-heat}
We  illustrate the use of the
Fourier theory by applying it to the heat equation on $X$.
Denote by $L_X$ the Laplace operator
on $X$. It is known that $e_{\lambda,b}$
is an eigenfunction for $L_X$, for all $\lambda$, $b$,
with the eigenvalues
$-(|\lambda |^2 + |\rho|^2)$ for $\lambda\in i\fa^*$.

The {\it heat equation} on $X$ refers to the Cauchy problem
\begin{equation}\label{eq-heat}
L_X u(x,t)=\partial_t u(x,t)\qquad \mathrm{and}\qquad
u(x,0)=f(x)
\end{equation}
where $f\in L^2(X)$. The problem has a unique solution,
which is easily
found by using the Fourier transform in Theorem \ref{th-He}.
It is given by
\begin{equation}\label{eq-heatsolution}
u(x,t)=\frac{1}{|W|}\int_{i\fa^*\times B}e^{-t(|\lambda |^2+|\rho |^2)}
\widehat{f}(b, \lambda )e_{\lambda ,b }(x )
\, d\mu (\lambda ,b )\, .
\end{equation}
We denote by $h_t$ the \textit{heat kernel} on $X$ defined by
$\hat h_t(\lambda,b)=e^{-t(|\lambda
|^2+|\rho |^2)}$ for all $\lambda,b$, or equivalently,
\begin{equation}\label{eq-heatkernel}
h_t(x)=\frac{1}{|W|}\int_{i\fa^*} e^{-t(|\lambda
|^2+|\rho |^2)}\varphi_{\lambda}(x)\,
d\mu (\lambda )\,
\end{equation}
(see \cite{Ga68}).

The convolution product of a function $f$ on $X$ with a
$K$-invariant function $h$ on $X$ (both viewed as functions
on $G$), is again a function on $X$.
Moreover, since
$\pi(f*h)=\pi(f)\pi(h)$, it follows that
$$(f*h)^\wedge(\lambda,b)=\hat f(\lambda,b)\hat h(\lambda)$$
It follows that  we can write, for all $f\in L^2(X)$,
\begin{equation}\label{eq-heatsol}
u(\cdot ,t)=f*h_t.
\end{equation}

We define the {\it heat transform}  as the map
$f\mapsto H_tf=f*h_t$ that associates the solution
at time $t>0$ to the initial function $f$. It follows from
the Fourier analysis that the heat transform is injective.
Notice that the map $H_t$ is a $G$-equivariant bounded operator
from the space $L^2(X)$ to itself. We equip the image $\Im (H_t)$,
not with the norm of $L^2(X)$, but with the norm that makes
the heat transform a unitary isomorphism.
As a consequence we obtain the following
\begin{observation}\label{th-ImHt}
Let $d\mu_t (\lambda )=e^{2t(|\lambda |^2+|\rho |^2)}\, d\mu (\lambda )$.
Then, as a representation of $G$, the image of the heat transform
decomposes as
\begin{equation}\label{eq-imageofHt}
(\Im (H_t),L)\simeq \int_{i\fa_+^*}^\oplus (H_{-\lambda } ,\pi_{-\lambda})\, d\mu_t(\lambda )\, .
\end{equation}
\end{observation}

Notice also the semigroup property $H_t(H_sf)=H_{s+t}f$
which follows from (\ref{eq-heatkernel}).

\section{The complex crown of $X$ and the space $\cHC$}\label{s-cHC}
\noindent
In this section we discuss some aspects of the interplay between
complex geometry and the harmonic analysis on
$X$. We introduce the crown, and we construct a $G$-invariant
Hilbert space of holomorphic functions on it. The
construction is motivated by the analysis in \cite{GKO04a} where
a similar construction was carried out on a subdomain
$\Cr_j\subseteq \Cr$. There the purpose was to obtain a Hardy space
realization of a part
of the most continuous spectrum of a pseudo-Riemannian
symmetric space $G/H_j$,
which is embedded in the boundary of $\Cr_j$.
In the present paper, our purpose is obtain a holomorphic
model that carries all the representations in
the Plancherel decomposition of $L^2(X)$.

The main references for Subsection \ref{ss-CROWN}  are
\cite{GK02,KS05,GKO04a}.  For a nice overview article
see \cite{F03a}.

\subsection{The complex convexity theorem}\label{ss-CROWN}
Let
$$\Omega=\{X\in \fa\mid (\forall \alpha \in \Delta )\, |\alpha (X)|<\pi/2\}\, .$$
For $\emptyset\not=\omega\subset\fa$ set
$T(\omega )=\fa\oplus i\omega$, $A (\omega )=\exp T(\omega)$ and
$X(\omega )=G A(\omega )\cdot x_o= G\exp i\omega \cdot x_o$.
Then $X(\omega )$ is $G$-invariant
by construction.
The set $\Cr:=X(\Omega )$, introduced in \cite{AG90}, is called the {\it complex crown}
of $X$
or the Akhiezer-Gindikin domain.
Note, that by Theorem 11.2 in \cite{He78}  it follows that
$\exp : T(\Omega ) \to A(\Omega )$ is a diffeomorphism.

The crown is an open $G$-invariant complex submanifold of $X_\C$
containing $X$ as a totally real submanifold.
Furthermore, the $G$-action
on $\Cr$ is proper,
$\Cr \subset N_\C A_\C\cdot x_o$ and $\Cr$ is Stein. The importance
of $\Cr$ for harmonic analysis comes from the fact,
that it allows a holomorphic extension of every eigenfunction of
the algebra of invariant differential operators.
The main step in proving this,
see \cite{KS05} Proposition 1.3, is
to show that each function $e_{\lambda, b}$ extends to a holomorphic
function on $\Cr$ and then to use the affirmative solution to
the Helgason conjecture \cite{KKMOOT78}.

The important complex convexity theorem of Gindikin-Kr\"otz \cite{GK02} is the inclusion
$\subseteq$ of the following theorem. The equality was recently
established by Kr\"otz-Otto in \cite{KOt06}.

\begin{theorem}\label{th-convexity}
Let $Y\in \Omega$, then
\begin{equation}\label{eq-convex2}
a (\exp (iY) G)= A\exp ( i\,\conv (W\cdot Y))
\, .
\end{equation}
\end{theorem}
Here $\conv$ stands for convex hull.
Note that (\ref{eq-convex2}) follows form Theorem 5.1
in \cite{KOt06} by $a(\exp iYg)=a(\exp iY k(g))a(g)$.

\subsection{The space $\cHC$}\label{ss-cHC}
We define the following function
$\omega: i\fa^*\times \fa \to \R^+$.
\begin{equation}\label{eq-omegaY}
\omega (i\nu,Y):=\frac{1}{|W|}\sum_{w\in W} e^{2\nu(wY)}\in \R^+.
\end{equation}
Furthermore, we define
$$\omega (\lambda)=\sup_{Y\in\Omega}\omega (\lambda ,Y)\, $$
for $\lambda\in i\fa^*$, and we
define a measure $\mu_\omega$ on $i\fa^*\times B$ by $d\mu_\omega
(\lambda,b) = \omega (\lambda )d\mu (\lambda ,b)$.

\begin{lemma}\label{estimate for e}
Let $Y\in \Omega$. There exists $T\in\fa^+$,
and for each $g\in G$ a constant $C>0$,
such that
\begin{equation}\label{eq-estimate}
\left|\frac{e_{i\nu , b}(g\exp iY \cdot x_o)}{\sqrt{\omega (i\nu )}}\right|
\le Ce^{-\nu (T)}
\end{equation}
for all $\nu\in\fa^*_+$ and $b\in B$.
The element $T$  can be chosen locally uniformly
with respect to $Y$, and
the constant $C$ can be chosen locally uniformly
with respect to $Y$ and $g$.
In particular
$(\lambda ,b) \mapsto e_{\lambda ,b }(z)/\omega (\lambda )$
is in $L^2(i\fa_+^*\times B,d\mu_\omega )$ for all $z\in\Cr$.
\end{lemma}
\begin{proof} Let $b=kM\in B$. Using Theorem \ref{th-convexity}
we can write
$$a(\exp (-iY)g^{-1}k)=a\exp iZ$$
where
$Z\in -\conv (W\cdot Y)\subset \Omega$ and $a\in A$ are unique and
depend continuously on
$Y$, $g$ and $b$. In particular,
if $Q\subset G\times \Omega$ is compact, then there
exists a constant $C=C_Q>0$ such that
$a^{-\rho}\le C$ for all $(g,Y)\in Q$ and $b\in B$.
Hence, for all $\lambda =i\nu \in i\fa^*_+$:
\begin{equation*}
|e_{\lambda , b}(g\exp iY \cdot x_o)|
=|(a\exp(iZ))^{-i\nu-\rho}|
\le C e^{\nu (Z)}.
\end{equation*}
Let $U\subset \Omega$ be a compact neighborhood of $Y$.
Since $-\conv (W\cdot U)$ is compact and contained in $\Omega$,
we can find an element $T\in \fa^+$
such that, $-\conv (W\cdot U)+T\subset \Omega$.
Thus by (\ref{eq-omegaY})
$$\frac{1}{|W|}e^{2\nu (Z+T)}\le
\frac{1}{|W|}\sum_{w\in W}e^{2\nu (w(Z+T))}
=\omega(i\nu,Z+T)
\le \omega (\lambda)$$
for all $Z\in -\conv (W\cdot U)$.
Notice that $T$ was chosen independently of $g$ and $b$.
Hence
$$\frac{1}{\omega (\lambda )}\le |W|e^{-2\nu (Z)}e^{-2\nu (T)}\, .$$
It follows that there exists a constant $C>0$  as claimed such that
(\ref{eq-estimate}) holds.

As $|c(i\nu  )|^{-2}$ has a polynomial growth on $\fa^*$ it follows
that $|c(i\nu )|^{-2}e^{-\nu (T)}$ is bounded. The last statement follows as $\fa_+^*\ni \nu \mapsto
e^{-\nu (T)}$ is integrable for all $T\in\fa^+$, see \cite{FK},
p. 10.
\end{proof}

Denote by $\cHC$ the space of holomorphic functions $F:\Cr\to \C$
such that $F|_X\in L^2(X)$ and
$$\|F\|_{\cHC}^2:=\frac{1}{|W|}\int_{i\fa^*\times B}
|\widehat{F|_{X}}(\lambda ,b)|^2
\, d\mu_\omega (\lambda ,b)<\infty\, .$$

\begin{theorem}\label{th-3.1} The
space $\cHC$ is a $G$-invariant Hilbert space. The
action of $G$ is unitary and
$$(\cHC,L) )\simeq \int_{i\fa_+^*}^\oplus (H_{-\lambda},\pi_{-\lambda })
\, d\mu_\omega (\lambda )\, .$$
Furthermore, the following holds:
\begin{enumerate}
\item Let $F\in\cHC$ and $f=F|_X$. Then
\begin{equation*}
F(z)=\int \widehat{f} (\lambda ,b )e_{\lambda ,b}(z)\, d\mu (\lambda ,b)
\end{equation*}
for all $z\in \Cr$.
\item For each $\varphi\in L^2(i\fa_+^*\times B,d\mu_\omega)$
the function defined by
\begin{equation*}
F(z)=\int \varphi (\lambda ,b )e_{\lambda ,b}(z)\, d\mu (\lambda ,b)
\end{equation*}
belongs to $\cHC$ and has $\widehat{F|_X}=\varphi$.
\item The point evaluation maps
$$\cHC\ni F\mapsto F(z)\in\C$$
are continuous for all $z\in \Cr$.
\item The
reproducing kernel of  $\cHC)$ is given by
\begin{eqnarray*}
K(z,w)&=&
\int_{i\fa_+^*\times B} \frac{ e_{\lambda ,b}(z)e_{-\lambda ,b}(\sigma (w)) }{\omega (\lambda )}\,
d\mu (\lambda ,b)\,\\
&=&\int_{i\fa_+^*} \frac{ \varphi_{\lambda}(\sigma (w)^{-1}z)}
{\omega (\lambda )}\,
d\mu (\lambda)\,
\end{eqnarray*}
where $\sigma$ is the conjugation introduced in Subsection \ref{ss-str}.
\end{enumerate}
\end{theorem}

\begin{proof}
We first establish the property (2). For $z\in\Cr$ let $f_z(\lambda ,b)=e_{\lambda ,b}(z)\omega (\lambda )^{-1}$.
Assume that $\varphi \in L^2(i\fa_+^*\times B,d\mu_\omega )$. Then
by Cauchy-Schwartz and the last part of Lemma \ref{estimate for e}
we have
\begin{eqnarray*}
\int_{i\fa^*_+\times B}|\varphi (\lambda ,b )
e_{\lambda ,b}(z)|\, d\mu (\lambda ,b) &=&\int_{i\fa^*_+\times B}|\varphi (\lambda ,b )
\frac{f_z (\lambda ,b(z)}{\omega (\lambda )}|\, d\mu_\omega (\lambda ,b)  \\
&\le &  \|\varphi \|_{L^2(d\mu_\omega  )}
\|f_z\|_{L^2(d\mu_\omega  )}\\
&<& \infty
\end{eqnarray*}
Recall that there is a compact neighborhood $U$ of $z$ such that $\|f_z\|$ is
uniformly bounded on $U$. It follows that the function
$$\Cr \ni z \mapsto G_\varphi(z):=\int \varphi (\lambda ,b )e_{\lambda ,b}(z)\, d\mu (\lambda ,b)\in \C$$
exists and is holomorphic.
Furthermore,
\begin{equation}\label{eq-pointevaluation}
|G_\varphi(z)|\le C\|\varphi\|_{L^2(d\mu_\omega )}\, .
\end{equation}
for all $\varphi\in  L^2(i\fa_+^*\times B,d\mu_\omega )$
with $C$ depending locally uniformly on $z$.
Now part (2) follows.

Let $F\in\cHC$ and $f=F|_X$.
Then it follows that
$$G_{\hat f}(z)=\int \widehat{f} (\lambda ,b )e_{\lambda ,b}(z)\, d\mu (\lambda ,b)$$
is holomorphic and satisfies $G|_X=F|_X$. Hence $G=F$.
This proves part (1), and
the part (3) now follows  from (\ref{eq-pointevaluation}).

It follows from part (3) that
there exists a reproducing kernel
$K(z,w)=K_w(z)$ for $\cHC$.
Let $F\in \cHC$. Then on the one hand
\begin{equation*}
F(w )=(F,K_w)_{\cHC}=
\int \widehat{F|_X}(\lambda ,b ) \overline{\widehat{K_w|_X} (\lambda ,b )}
\omega (\lambda )\, d\mu (\lambda ,b),
\end{equation*}
and on the other, using part (1)
$$F(w)=\int \widehat{F|_{X}}(\lambda ,b ) e_{\lambda, b}(w) \, d\mu (\lambda ,b).$$
It follows from part (2) that $\widehat{F|_{X}}(\lambda ,b )$
can be any function in $L^2(i\fa_+^*\times B,d\mu_\omega)$.
Since  $e_{\lambda, b}/\omega$ also belongs to this space
by Lemma \ref{estimate for e},
it follows that
$$\widehat{K_w|_X} (\lambda ,b )=\overline{ e_{\lambda, b}(w)}/ \omega (\lambda ) =
e_{-\lambda ,b}(\sigma (w))/\omega (\lambda ).$$
Thus
\begin{eqnarray*}
K(z,w)&=&(K_w,K_z)_{\cHC}\\
&=&\int \frac{e_{-\lambda ,b}(\sigma (w)) e_{\lambda ,b}(z)}{\omega (\lambda )^2}\, \omega (\lambda ) d\mu\\
&=&\int \frac{ e_{\lambda ,b}(z)e_{-\lambda ,b}(\sigma (w))}{\omega (\lambda )}\,  d\mu\,.
\end{eqnarray*}
It suffices to establish the final formula in (4) for $z,w\in X$.
Moreover, by $G$-invariance of the kernel we may assume $w=e$.
Then the formula follows immediately from (\ref{spherical function})
(see also Theorem 1.1 in \cite{He94}, p. 224).
\end{proof}

\begin{remark} It is clear that the space $\cHC^K$ also is a reproducing kernel Hilbert space. The reproducing
kernel is obtained by averaging over $K$,
$$K_{\cHC^K}(z,w)=\int_{i\fa_+^*} \frac{ \varphi_{\lambda}(z)
\varphi_{-\lambda}(\sigma (w)) }{\omega (\lambda )}\,
d\mu (\lambda)\, ,
$$
as stated in \cite{GKO06a}, Proposition 8.7.
\end{remark}

\section{The image of the Segal-Bargmann Transform}\label{s-imageHt}
\noindent
In this section we introduce the \textit{Segal-Bargmann transform} on $X$
and give two characterizations of its image as a $G$-invariant
Hilbert space of holomorphic functions on $\Cr$, both
different from the one given in \cite{KOS05}.
The first characterization is a natural extension of Observation \ref{eq-imageofHt}.
The second characterization
uses the normalized Radon transform $\Lambda$ from \cite{He70}.
A similar result was proved in \cite{OS06} for the
Segal-Bargmann transform related
to positive multiplicity functions, but only for the
$K$-spherical case.

In the last part of the section we
show that the normalized Radon transform allows an analytic extension
as a unitary
isomorphism of $\cHC$ into a function space over a domain in
the complexified horocycle space
$\Xi_\C$.
A similar
construction for the Hardy space on a subdomain $\Cr_j$ was given
in \cite{GKO06a}.

\subsection{The Segal-Bargmann transform on $X$}
The description of the heat transform in Subsection \ref{ss-heat}
implies that the heat kernel $h_t$, as well as every function $H_tf$,
$f\in L^2(X)$, extends to a holomorphic function on
$\Cr$ (see also \cite{KS05}, Prop. 6.1).
In fact, it follows from Theorem \ref{th-3.1} (2),
that these extensions belong to $\cHC$. We shall denote the
holomorphic extensions by the same symbols.
The map $$H_t : L^2(X)\to\cHC \subset\cO (\Cr )$$ is the
\textit{Segal-Bargmann transform}. Note, $\|H_tf\|_{\cHC}\le
\|f\|_{L^2}$.

Our first description of the image of $H_t$ in $\cO (\Cr )$ is given
by the following.
\begin{theorem}\label{th-ImageHt2}
The image $H_t(L^2(X))$ of the Segal-Bargmann transform
is the space $\cO_t(\Cr )$ of
holomorphic functions $F$ on $\Cr$ such that $F|_X\in L^2(X)$ and
$$\|F\|_t^2:=\int_{i\fa_+^*\times B} |\widehat{F|_{X}}(\lambda ,b)|^2\, d\mu_t(\lambda ,b)<\infty\, .$$
Here $d\mu_t$ is the measure defined in
Observation \ref{th-ImHt}.

Furthermore,
point evaluation is continuous on $\cO_t(\Cr )$, and the
reproducing kernel  is given
by
$$K_t(z,w)=h_{2t}(\sigma(w)^{-1}z)\, .$$
\end{theorem}

\begin{proof} The first statement
follows from Observation \ref{th-ImHt}
and Theorem \ref{th-3.1}. The  reproducing kernel is obtained by
a standard argument, using the semigroup property of the
convolution with $h_t$. Let $F=f*h_t$ and $z \in\Cr$. Note
that $L_zh_t$ is well defined for $z\in \Cr$
as $h_t$ is $K$-biinvariant. We have
\begin{eqnarray*}
F(w)&=&f*h_t(w)\\
&=&(f,L_{\sigma (w )}h_t)_{L^2(X)}\\
&=&(H_tf,H_t(L_{\sigma (w )}h_t))_{\cO_t}\, .
\end{eqnarray*}
Hence point evaluation is continuous and
$$K_w(z)=H_t(L_{\sigma (w )}h_t)(z)=(L_{\sigma (w)}h_t)*h_t(z)=h_{2t}(\sigma (w)^{-1}z)\, .$$
\end{proof}

\subsection{The Radon transform on X and the Segal-Bargmann transform}\label{ss-Radon}
Let $\Xi=G/MN$ and put
$\xi_o =eMN\in \Xi$.  By the Iwasawa decomposition it follows that
\begin{equation}\label{eq-polarXi}
B \times A\simeq \Xi\, \quad (kM,a)\mapsto ka\cdot \xi_o
\end{equation}
is a diffeomorphism. A subset $\xi \subset X$ is said to be a
\textit{horocycle} if there exists $x\in X$ and $g\in G$ such
that
$\xi =gNx$. It is well known, see \cite{He63,He70}, that
the map $ g\xi_0\mapsto gNx_0$ is a bijection
of $\Xi$ onto the set of horocycles. Using this identification
the space of horocycles becomes an analytic manifold with a
transitive $G$ action.

The Radon transform $\cR (f)(g\cdot \xi_o)
=\int_N f(gn\cdot x_o)\, dn$ is a $G$-intertwining
operator $C_c^\infty (X)\to C^\infty (\Xi )$.
The importance of this observation comes
from the fact that the regular representation of $G$
on $L^2(\Xi)$ is much easier to decompose than that on
$L^2(X)$. As induction commutes with direct integral,
induction by stages shows that
$(L^2(\Xi ) ,L)\simeq \int_{i\fa^*}^\oplus
(H_\lambda ,\pi_\lambda )\, d\lambda $,
see \cite{M76}, p. 284 and 287.
In fact, let $\chi_\lambda (man )=a^{\lambda}$.
Denote by $\epsilon$ the trivial representation
of $MN$ and by $L_\Xi$ the regular representation
of $G$ on $L^2(\Xi )$. As $MAN/MN\simeq A$ and  $MN$ acts trivially on
$L^2(MAN/MN)\simeq L^2(A)$, it follows that
\begin{eqnarray*}
L_\Xi&\simeq&\ind_{MN}^G\epsilon\\
&\simeq&\ind_{MAN}^G\ind_{MN}^{MAN}\epsilon\\
&\simeq&\ind_{MAN}^G\int_{i\fa^*}^\oplus \chi_\lambda\, d\lambda\\
&\simeq&\int_{i\fa^*}^\oplus \pi_\lambda\, d\lambda\, .
\end{eqnarray*}

Equation (\ref{eq-HaarOnG}) implies that (up to a constant) the $G$-invariant measure on
$\Xi$ is given by
$$\int_\Xi f(\xi)\, d\xi = \int_{B}\int_A f(ka\cdot \xi_o)\, a^{2\rho}\, dadk\, .
$$
It follows that
\begin{equation}\label{Xi to B times A}
L^2(\Xi)\ni f\mapsto [(kM,a)\mapsto a^\rho f(ka\cdot\xi_0)]\in L^2(B\times A)
\end{equation}
is a unitary isomorphism. We also note that by
(\ref{eq-HaarOnG}),
\begin{equation}\label{eq-FourierXandA1}
\hat{f}(\lambda ,kM)=\int_{AN}f(kan\cdot x_o)a^{-\lambda + \rho}\, dnda
=\int_A a^\rho \cR(f)(kM,a)a^{-\lambda }\, da
\end{equation}
for $f\in C_c(X)$.
It is therefore natural to define a $\rho$-twisted Radon transform  by
\begin{equation}\label{eq-defRrho}
\cR_\rho (f)(b,a)=a^\rho \cR(f)(b,a)\, .
\end{equation}
Notice that this is then an intertwining operator for the regular action
of $G$ on functions over $X$ and $\Xi$, respectively, when the action
on functions over $\Xi$ is transferred to a $\rho$-twisted action
on functions over $B\times A$ through
(\ref{Xi to B times A}), that is,
\begin{equation}\label{rho twisted action}
(g\cdot\phi)(b,a):=a(g^{-1}b)^{-\rho}\phi(k(g^{-1}b),a(g^{-1}b)a).
\end{equation}

Identifying $L^2(B\times A)$ with $L^2(A,L^2(B))$ in a natural way,
equation (\ref{eq-FourierXandA1}) now reads
\begin{equation}\label{eq-FourierXandA2}
\hat{f}_\lambda =\cF_A (\cR_\rho (f))(\lambda)\, .
\end{equation}
Note, if $f$ is $K$-invariant the $\rho$-twisted Radon transform reduces to
the Abel transform $f\mapsto F_f$ introduced in
\cite{HC58}, p.~261, and conjectured  to be injective.
The proof of that conjecture was the final step towards the
Plancherel formula, see \cite{HC66}, p.~4. It follows from (\ref{eq-FourierXandA2})
and Theorem \ref{th-He} that $\cR_\rho$ is injective also without
the assumption of $K$-invariance.

Denote by $\Psi_\fa$ the multiplication operator $F\mapsto \frac{1}{c(-\cdot )}F$
and define a pseudo-differential operator $\Psi$ on $A$ by
$$\Psi =\cF_A^{-1}\circ \Psi_\fa\circ \cF_A\, .$$
Recall that by the
Gindikin-Karpelevich product formula
for the $c$-function, it
follows that  $\Psi$ is a differential operator if and only if all
multiplicities $m_\alpha$ are even.

Let $L^2_W (B\times i\fa^*, |W|^{-1}db\,d\lambda)
\simeq L^2_W (i\fa^*, L^2(B); |W|^{-1}\, d\lambda)$ be the
space of functions in  $L^2 (B\times i\fa^*, |W|^{-1}db\,d\lambda)$
satisfying the intertwining
relations
\begin{equation}\label{eq-WrelationLambda}
c(-w\lambda )F(\cdot , w\lambda )=c(-\lambda ) \cA (w,-\lambda ) F (\cdot ,\lambda )\, ,
\end{equation}
for all $w\in W$,
cf. Theorem \ref{th-He}. Note, if $F$ is $K$-invariant, then
(\ref{eq-WrelationLambda})
amounts to
$c (-w\lambda )F(w\lambda )=c(-\lambda ) F(\lambda )$ for all $w\in W$.
Let $d\tau (b,a)=|W|^{-1}db\,da$ on $B\times A$.
The preimage of $L^2_W (B\times i\fa^*, |W|^{-1}db\,d\lambda)$
under the Fourier
transform
of $A$ is denoted by  $L^2_W(B\times A,d\tau )$. It is often
helpful to identify this space with $L^2_W (A,L^2(B);|W|^{-1}\, da) $.
Several formulas like the definition of the Hardy space $\cHXi$
in (\ref{def-hardy2}) below are then analogous to the Euclidean case,
cf. \cite{SW71}, Chapter III, Section 2.

By Theorem \ref{th-He} and (\ref{eq-FourierXandA2}) we have the
following commutative diagram
\begin{equation}\label{def-Lambda}
\xymatrix{L^2(X )\ar[d]_{\mathcal{F}}&
 \ar[r]^{(\mathrm{id}\times \Psi )\circ \cR_\rho}&
&
L^2_W(B\times A,d\tau )\ar[d]^{\mathrm{id}\times \mathcal{F}_A}
\\
L^2_W(B \times \fa^*,|W|^{-1}d\mu )& \ar[r]_{\mathrm{id}\times \Psi_\fa}&
& L^2_W(B\times \fa^* ,|W|^{-1}db\,d\lambda )
}\,
\end{equation}
The vertical maps and the lower horizontal map are unitary isomorphisms.
It follows that that the linear operator
\begin{equation}\label{defLambda}
\Lambda :=(\mathrm{id}\times \Psi)\circ \cR_\rho : L^2(X) \to L^2_W (B\times A,d\tau )
\end{equation}
is an unitary isomorphism. By Lemma 3.3 in \cite{He70}, p. 42, we know that
$\Lambda$ is an intertwining
operator.  Here the action of $G$ on
$L^2(B\times A,d\tau )$ is the $\rho$-twisted action
defined through the identification
of $B\times A$ with $\Xi$ (see (\ref{rho twisted action})).

By Theorem 5.3 in
\cite{He63} it follows that
for each  invariant differential operator $D$ on $X$,
there exists a differential operator $\widetilde{D}$ on $A$ such
that $\Lambda (Df)(b,a)=\widetilde{D}_a\Lambda (f)(b,a)$.
Here the subscript indicates that
$\tilde{D}$ acts on the $a$ variable only. In particular,
this applies with $D$ equal to the Laplace operator $L_X$.
Tracing the commutative diagram (\ref{defLambda})
and using that $\widehat{L_Xf}(\cdot , \lambda  )=-(|\lambda |^2+|\rho |^2)
\hat{f}(\cdot ,\lambda  )$,
$\lambda \in i\fa^*$,
it is easily seen that $\widetilde{L_X}=L_A-|\rho |^2$, i.e., for $f$ sufficiently smooth
\begin{equation}\label{eq-LaplaceRrho}
\Lambda (L_Xf)=(L_A-|\rho |^2)\Lambda (f)\, .
\end{equation}

Let $r=\dim A$.

\begin{lemma}\label{le-heatOnA} Let $f\in L^2(X)$. Then
$e^{t|\rho |^2}\Lambda (H_t f)$  solves the heat equation
on $A$ with initial value $\Lambda (f)\in L^2(B)$. In particular, the map
$\fa\ni X\mapsto \Lambda (H_t f)(\cdot ,\exp X)\in L^2(B)$ extends to a holomorphic function
on $\fa_\C$, again denoted by $\Lambda (H_t f)$ such that
$$|W|^{-1}(2\pi t)^{-r/2}
\int_{B\times \fa_\C}|e^{t|\rho|^2}\Lambda (H_t f)(b,\exp (X+iY))|^2\, e^{-|Y|^2/2t} \, db
dXdY<\infty\, .$$
\end{lemma}

\begin{proof} See the proof of Lemma 2.5 in \cite{OS06}.
\end{proof}
Let $\cF_{W,t}(B\times \fa_\C )$ denote the space of holomorphic functions on
$B\times \fa_\C$ such that
$$\fa_\C\ni Z\mapsto F(\cdot ,Z)\in L^2(B)$$
is holomorphic with $F\circ (\id,\log )\in L^2_W (B\times A)$, and satisfies
\begin{eqnarray*}
\|F\|^2_t&=& |W|^{-1}(2\pi t)^{-r/2}
\int_{\fa_\C^*}\int_B |F(b, X+iY )|^2\, e^{-|Y|^2/2t}\, dX dY\\
&=&|W|^{-1}(2\pi t)^{-r/2}\int_{\fa_\C^*}\|F(X+iY )\|^2_{L^2(B)}\, e^{-|Y|^2/2t}\, dXdY <\infty\, .
\end{eqnarray*}
Thus, $\cF_{W,t}(B\times\fa_\C)$ is analog to a $L^2(B)$ valued Fock space describing the image
of the Segal-Bargmann transform on the Euclidean space $\fa$ with the addition of
the Weyl group relations derived from (\ref{eq-WrelationLambda}). For $t>0$ define
$\Lambda_t : \cO_t(\Cr )\to \cF_{W,t}(B\times\fa_\C)$ in the following way.
Let $F\in \cO_t(\Cr )$. By Theorem \ref{th-ImageHt2} there exists a
unique $f\in L^2(X)$ such that $F|_X=H_tf$. Let $\Lambda_t(F)$ be the holomorphic extension of
$e^{t|\rho |^2}\Lambda (H_tf)$ given by Lemma \ref{le-heatOnA}. By the
same lemma the Weyl group relations are satisfied and  $\|\Lambda_t(F)\|_t<\infty$.
Thus $\Lambda_t(F)\in\cF_{W,t}(B\times\fa_\C)$.
The following theorem gives an alternative
description of the Hilbert space $\cO_t(\Cr)$.

\begin{theorem}
\label{th-image2} The map $\Lambda_t  :\cO_t (\Cr )\to \cF_{W,t}(B\times \fa_\C) $
is an unitary isomorphism. Furthermore, let $F\in \cO_t(\Cr )$. Define
$f\in L^2(X)$ by applying $\Lambda^*$ to the function on $B\times A$ given by
$$(b,a)\mapsto (4\pi t)^{-r/2}\lim_{R\to \infty}\int_{|Y|\le R}\Lambda_t(F)(b,\log a + iY)
e^{-|Y|^2/4t}\, dY\, .$$
Then $H_t(f)=F$.
\end{theorem}

\begin{proof} The proof is a simple adaption of the standard argument for
$\R^r$, as described in \cite{H04a},  to the $L^2(B)$-valued case.
\end{proof}

\subsection{Holomorphic properties of the normalized Radon transform}\label{ss-heat2}
\noindent
Consider as before $X$ as a subset of $X_\C=G_\C/ K_\C$.
A  complex horocycle
in $X_\C$ is a set of the form $gN_\C x_o\subset X_\C$
for some $g\in G_\C$, see \cite{GKO06a}. Let $\xi_{o}^\C=N_\C\cdot x_o$, then
$$\Xi_\C =\{g \cdot \xi_{o}^\C \subset G_\C/K_\C\mid g\in G_\C\}\simeq G_\C/M_\C N_\C$$
is the set of complex horocycles.
The map
$$\Xi \ni g\cdot \xi_{o}\mapsto g\cdot \xi_{o}^\C\in \Xi_\C
,\quad g\in G$$
is well defined and injective.

Define
$$\Xi (\Omega )=G\exp  i\Omega\cdot \xi_o=KA(\Omega )\cdot \xi_o\subset \Xi_\C\, .$$

Then $ \Xi (\Omega )\simeq B\times A(\Omega) $
is a $G$-invariant CR-submanifold of $\Xi_\C$.
Let $\cHXi$ be the space of function $F : \Xi (\Omega )\to \C$ such that
the map $A (\Omega )\ni z\mapsto F(\cdot ,z)=F_z\in L^2(B)$
is holomorphic and for each fixed $Y\in \Omega$ the function $(b,a)\mapsto
F (b, a\exp iY)$ is
in $L^2_W (B\times A,d\tau)$ and
\begin{eqnarray}\label{def-hardy2}
\|F  \|_{\cHXi }^2&:=& \sup_{Y\in \Omega }|W|^{-1} \int_{B\times A}
|F(b,a\exp iY)|^2\, dbda
\\
&=&
|W|^{-1} \sup_{Y\in \Omega } \int_{A} \|F_{a\exp iY} \|^2_{L^2(B)}\, da<\infty\, .\nonumber
\end{eqnarray}
We define an action of $G$ on $\cHXi$ by holomorphic extension of the
$\rho$-twisted action in (\ref{rho twisted action}),
that is
\begin{equation}
(g\cdot F)(ka\cdot\xi_0^\C):=a(g^{-1}k)^{-\rho}F(g^{-1}ka\cdot\xi^\C_0)
\end{equation}
for $g\in G$, $k\in K$ and $a\in A(\Omega)$. Notice that
$(g\cdot F)|_{B\times A}=g\cdot(F|_{B\times A})$.

\begin{lemma}\label{le-HCHXI1} Let $F=H_t(f)\in \cHC$ where
$f\in L^2(X)$, $t>0$. Then $a\mapsto \Lambda (F|_X)(\cdot ,a)$ extends to
a holomorphic $L^2(B)$ valued function on $A(\Omega )$, also denoted
by $\Lambda (F|_X)$, which belongs to $\cHXi$ and satisfies
$$\|\Lambda (F|_X )\|_{\cHXi}=\|F\|_{\cHC}\, .$$
Moreover, the map $f\mapsto \Lambda (F|_X)$ is intertwining for the actions of
$G$.
\end{lemma}
\begin{proof}
Let $\varphi = F|_X$. It follows by Lemma \ref{le-heatOnA}
that $a\mapsto \Lambda (\varphi )(\cdot ,a)$
extends to a holomorphic $L^2(B)$-valued function
on $A(\Omega )$. In fact
$$\Lambda (\varphi )(\cdot , \exp (X+ iY))=e^{-t|\rho |^2}\Lambda (f)*_A h_t^A (\cdot ,\exp (X+ iY) )$$
where $h_t^A(\exp X )=(4\pi t)^{-r/2}e^{-|X |^2/4t}$ is the heat
kernel on $A$
and the convolution is on the abelian group $A$.
For $Y\in \Omega$ the function
$g_Y: (b,a)\mapsto \Lambda (\varphi )(b ,a\exp iY)$ is
in $L^2_W (B\times A,d\tau )$.
By the explicit formula for $h_t^A$ there exists a
positive constant $C>0$ such that for $a\in A$ and $Y\in\Omega$
\begin{equation}\label{eq-estimatForLambdaVarphi}
\|g_Y(\cdot ,\exp X )\|_{L^2(B)}\le C e^{-(|X |-1)^2/4t}e^{|Y|^2/4t}
\le C_1e^{-(|X|-1)^2/4t}
\end{equation}
where
$$C_1=C\sup_{Y\in\overline{\Omega}}e^{|Y|^2 /4t}\, .$$

Let $g_{b,Y} (a)=g_Y(b,a)$. The estimate (\ref{eq-estimatForLambdaVarphi})
allows us to change the path of integration to derive
$$\cF_A (g_{b,Y})(\lambda) = \cF_A (g_{b,0}) (\lambda )e^{i\lambda (Y)}\, .$$
Thus the integral over $B\times A$ in (\ref{def-hardy2}) is
\begin{eqnarray}
\int_{B}\!\!\int_A |g_Y(b,a )|^2\, da db& =&\int_B\!\! \int_{i\fa^*}
|\cF_A(g_{b,Y} ) (\lambda )|^2\, \, d\lambda db\nonumber\\
&=&  \int_{i\fa^*}\!\!\int_B |\cF_A( g_{b,0}) (\lambda )|^2\,
e^{2i\lambda (Y)}\, db d\lambda\nonumber\\
&=& \int_{i\fa^*}\!\!\int_B |\cF_A (\Lambda (\varphi )) (b,\lambda )|^2\, db\, e^{2i\lambda (X)}\, d\lambda\nonumber\\
&=&\int_{B\times i\fa^*} |\widehat{\varphi }(b,\lambda )|^2e^{2i\lambda (Y)}\, \frac{dbd\lambda}{|c (\lambda )|^2}
\label{eq-normeq}
\end{eqnarray}
where (\ref{eq-normeq}) follows from the definition of $\Lambda$ in (\ref{def-Lambda}).
According to (\ref{eq-IntertwRel}) we have $\cA (w,-\lambda )
\widehat{\varphi }(\cdot ,\lambda )
= \widehat{\varphi }(\cdot, w\lambda )$.
Hence $\int_B |\widehat{\varphi }(b,\lambda )|^2\, db$ is $W$-invariant
as the intertwining operator $\cA (w,-\lambda )$ is unitary. Summing
over the Weyl group and
using that $|c(\lambda )|^{-2}$ is $W$-invariant,
we obtain
$$\int_{B\times i\fa^*} |\widehat{\varphi }(b,\lambda )|^2
e^{2i\lambda (Y)}\, \frac{dbd\lambda}{|c (\lambda )|^2}
=\int_{B\times i\fa^*} |\widehat{\varphi }(b,\lambda )|^2\,
\omega (\lambda  ,-Y)\,
d\mu (b,\lambda) \, .
$$
Divide by $|W|$ and take  the supremum over $Y\in\Omega=-\Omega$ to get that
the norms are equal.

The intertwining property of the map follows from the
corresponding properties for $H_t$ and for $\Lambda$ on $B\times A$.
The latter property was remarked below
(\ref{defLambda}).
\end{proof}

Let $F\in \cHC$ and $\varphi =F|_X$.
Let  $t_n\to 0$, $t_n>0$, and view $\varphi_n:=H_{t_n}\varphi$ as an
element of $\cHC$. Then
$$
\|\varphi_n-F\|^2_{\cHC}=\int_{B\times i\fa_+^*}|e^{-t_n(|\lambda |^2+|\rho |^2)}-1|^2|\hat{\varphi}(b,\lambda )|^2
\omega (\lambda )\, d\mu (b,\lambda)\, .$$
As $(b,\lambda )\mapsto |\hat{\varphi}(b,\lambda )|^2
\omega (\lambda )$ is integrable with respect to $d\mu $ it follows by
the
Lebesque dominant convergence theorem that $\lim_{n}\varphi_n\to F$ in
$\cHC$. By Lemma \ref{le-HCHXI1} it follows that $\lim_{n}\Lambda (H_{t_n}\varphi )$
exists in $\cHXi$ and is independent of the sequence $t_n$. Define $\tilde{\Lambda }: \cHC\to \cHXi$ by
$$\tilde{\Lambda}(F)=\lim_{n\to \infty}\Lambda (H_{t_n}\varphi )\, .$$
\begin{theorem}\label{th-HCIsoHXI} The map $\tilde{\Lambda }:\cHC\to
  \cHXi$ is an unitary intertwining
isomorphism.
\end{theorem}
\begin{proof} We only have to show that $\tilde{\Lambda}$ is surjective. Let
$F\in \cHXi$. Then $F|_{B\times A}\in L^2_W(B\times A,d\tau )$. Define
$f=\Lambda^*(F|_{B\times A})$. Then the argument in the proof of Lemma
\ref{le-HCHXI1}
shows that
$$|W|^{-1}\int_{B\times i\fa^*}|\hat{f}(b,\lambda )|^2\omega (\lambda )\, d\mu (b,\lambda )
=\|F\|^2_{\cHXi}<\infty\, .$$
In particular it follows that $f$ extends to a holomorphic function on $\Cr$, denoted by $\varphi$, and
$\varphi\in \cHC$. By construction we have $\Lambda (\varphi|_X)= F|_{B\times i\fa^*}$. By
construction it follows now easily that $\tilde{\Lambda} (\varphi )=F$.
\end{proof}

\def\paper{}\def\journal{}

\end{document}